\documentclass{article}

\usepackage{amssymb,latexsym,amsmath}

\usepackage[pdftex]{graphicx}
\usepackage{float}

\usepackage{hyperref, xcolor}

\hypersetup{backref, colorlinks=true, 
linkcolor=red, 
linkbordercolor = red,          
citecolor = blue}

\DeclareGraphicsExtensions{.jpg}

\begin{document}

\textwidth16.5cm

\voffset-2.5cm

\textheight20.8cm

\pagestyle{plain}

\newtheorem{theorem}{Theorem}[section]

\newtheorem{proposition}[theorem]{Proposition}

\newtheorem{lemma}[theorem]{Lemma}

\newtheorem{corollary}[theorem]{Corollary}

\newtheorem{definition}[theorem]{Definition}

\newtheorem{remark}[theorem]{Remark}

\newtheorem{exempl}{Example}[section]

\newenvironment{exemplu}{\begin{exempl}  \em}{\hfill $\square$

\end{exempl}}

\renewcommand{\contentsname}{ }

\title{Zipper logic}

\author{Marius Buliga \\
    Institute of Mathematics of the Romanian Academy \\ P.O. BOX 1-764, RO 014700, Bucharest, Romania \\ 
\texttt{Marius.Buliga@gmail.com}}

\date{20.05.2013}

\maketitle

\begin{abstract}
Zipper logic is a graph rewrite system, consisting in only local rewrites on a class of zipper graphs. Connections with the chemlambda artificial chemistry and with knot diagrammatics based computation are explored in the article. 
\end{abstract}

\section{Introduction}

Zipper logic is a graph rewrite system.  It consists in a class of graphs, called zipper graphs and a collection of moves (graph rewrites) acting on zipper graphs, introduced in section \ref{szippers}. Discover more  at the \href{http://chorasimilarity.wordpress.com/tag/zipper-logic/}{zipper logic tag} from the author's open notebook. 

Zipper logic is a variant of the chemlambda \cite{buligachem} \cite{alife} artificial chemistry, as explained in section \ref{schemlambda}. In section \ref{sturing} we prove that zipper logic is Turing universal, by showing that it can be used to implement the SKI combinatory logic. Because zipper logic, as chemlambda, has only local graph rewrites, we treat with special attention the processes of birth and death of zipper combinator graphs, which replace the  global fan-out rewrite move which is needed in (this graphical version of) combinatory logic. 

In section \ref{sknots} we turn to the exploration of relations between zipper logic and knot diagrammatics.

\section{Half-zippers, zippers and moves}
\label{szippers}

Let's start by defining the zipper graphs. Such a graph is made by the basic ingredients described in Fig.~\ref{define_nodes_1} and Fig.~\ref{4_chem}.

\begin{definition}
A zipper graph is an oriented graph which has as nodes: 
\begin{enumerate}
\item[-] $(-n)$ half-zippers, depicted in the first row of the Fig.~\ref{define_nodes_1}; for any natural number $n \geq 1$, a $(-n)$ half-zipper is a node with $n+2$ arrows, which are ordered, for convenience by numbering them with $0$, $0'$, $1$, ..., $n$, such that the arrow numbered by $0$ points to the node and the arrows numbered by $0'$, $1$, ..., $n$ point away from the node; 
\item[-] $(+n)$ half-zippers, depicted in the second row of the Fig.~\ref{define_nodes_1}; for any natural number $n \geq 1$, a $(+n)$ half-zipper is a node with $n+2$ arrows, which are ordered  by numbering them with $0$, $0'$, $1'$, ..., $n'$, such that the arrows numbered by $0$, $1'$, ..., $n'$ point to the node and the arrow numbered by $0'$ points away from the node; 
\item[-] $(n)$ zippers, depicted in the third row of the Fig.~\ref{define_nodes_1}; for any natural number $n \geq 1$, a $(n)$ zipper is a node with $2n+2$ arrows, which are separated into two disjoint sets of $n+1$ arrows;  the first set is formed by  the arrows numbered by $0$, $1'$, ..., $n'$, which  point to the node and the second set is formed by  the arrows numbered by $0'$, $1$, ..., $n$, which  point  away from the node; 

\begin{figure}[H]
\includegraphics[width=  105mm]{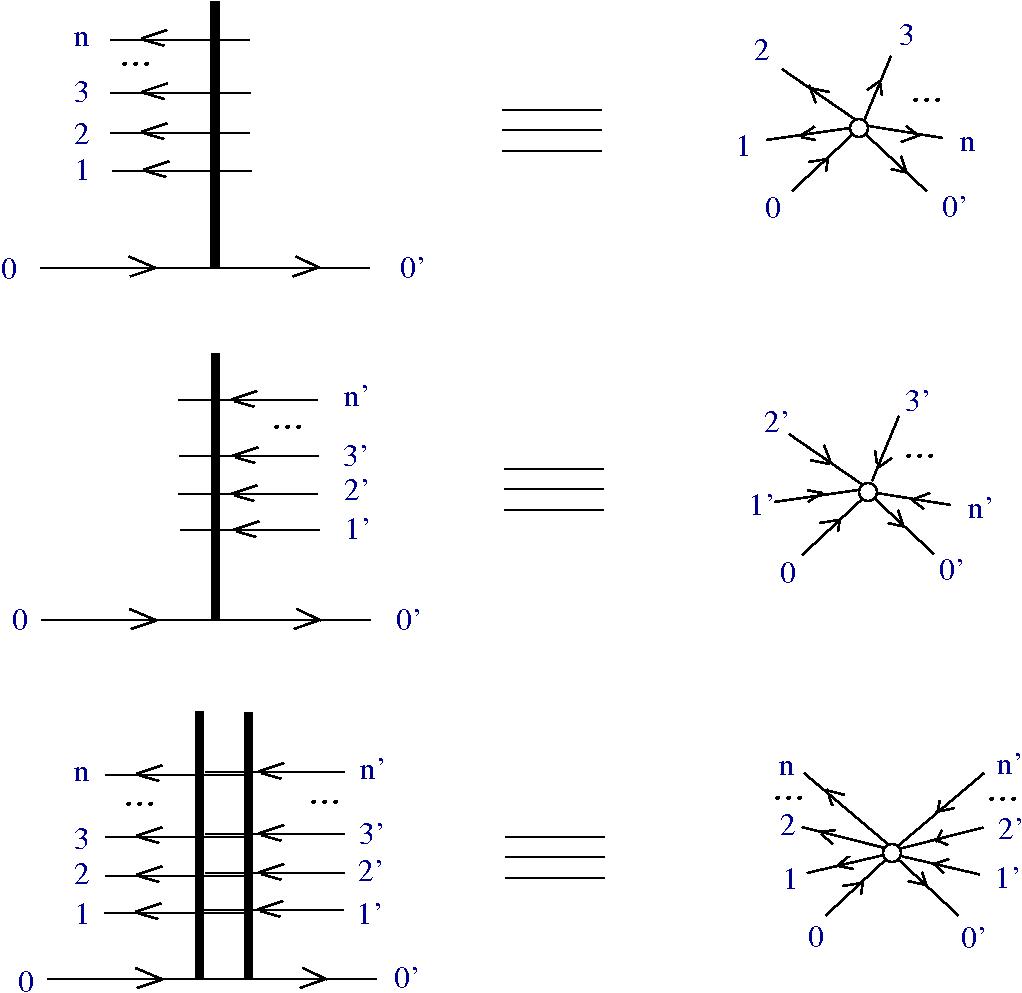}
\caption{Zipper nodes: (1st row) $(-n)$ half-zipper, (2nd row) $(+n)$ half-zipper, (3rd row) $(n)$ zipper.}
\label{define_nodes_1}
\end{figure}

\item[-] fanout and  fanin nodes, described in Fig.~\ref{4_chem} (a), (b), are trivalent nodes with a cyclic order of the arrows, 
\item[-] termination is a univalent node, with an arrow which points to the node, described  in Fig.~\ref{4_chem} (c).
\end{enumerate}
The arrows connect the nodes. We may have moreover arrows which point to one of the nodes but which have the origin not connected to any node, or we may have arrows with the origin connected to a node, but with the end free. Finally we may have arrows with both ends free, or loops, as described in Fig.~\ref{4_chem} (c). 

A zipper graph is made by a finite number of nodes, arrows and loops; it does not have to be connected. 

\label{defzippergraph}
\end{definition}

\begin{figure}[H]
      
     \includegraphics[width=  75mm]{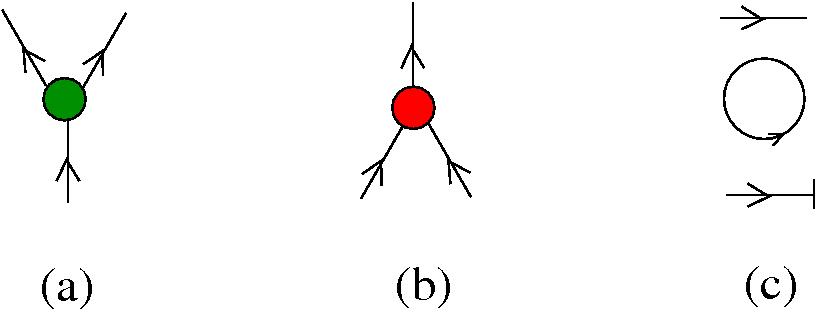}
     \caption{The other ingredients: (a) fanout node, (b) fanin node, (c) arrow, loop and termination node.}
     \label{4_chem}
 
\end{figure}

The graph rewrites of zipper logic are described in the following figures.

\begin{definition}
The moves of zipper logic are all reversible. They act on a bounded number of nodes (but due to the fact that they act on half-zippers or zippers with an unbounded number of arrows, they may act on any number of arrows). The figures which describe the moves contain only the region, or pattern, which is subjected to the respective move. They come in several families: 
\begin{enumerate}
\item[-] the CLICK moves, Fig.~\ref{click} transform pairs of half-zippers into a zipper and possibly a half-zipper; in the Fig.~\ref{click} is described a CLICK move which involves a $(-n)$ half-zipper and a $(+m)$ half-zipper with $m > n$; the other cases, namely $m=n$ and $m<n$ are not shown, but they are straightforward to imagine; 

\begin{figure}[H]
      
     \includegraphics[width=  105mm]{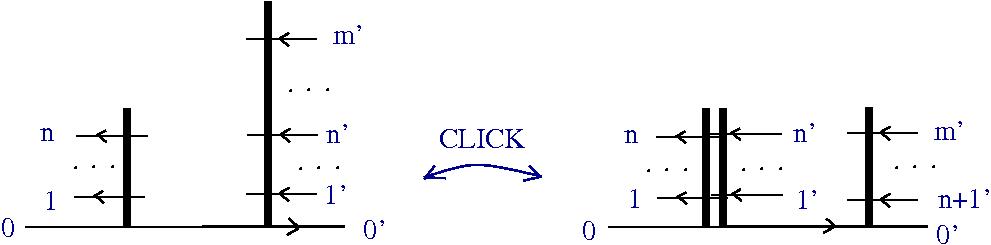}
     \caption{The CLICK move for $m > n$.}
     \label{click}
 
\end{figure}

\item[-] the ZIP move, Fig.~\ref{zip}, is the one which gives the name to the zipper logic, 

\begin{figure}[H]
      
     \includegraphics[width=  105mm]{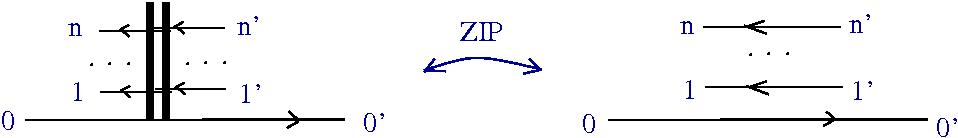}
     \caption{The ZIP move.}
     \label{zip}
 
\end{figure}

because it may be imagine as the act of unzipping a zipper, when seen from left to right, or to zip it, when seen from right to left; 

\item[-] the TOWER moves, Fig.~\ref{zipper_not_2}, serve to stack half-zippers, 

\begin{figure}[H]
      
     \includegraphics[width=  105mm]{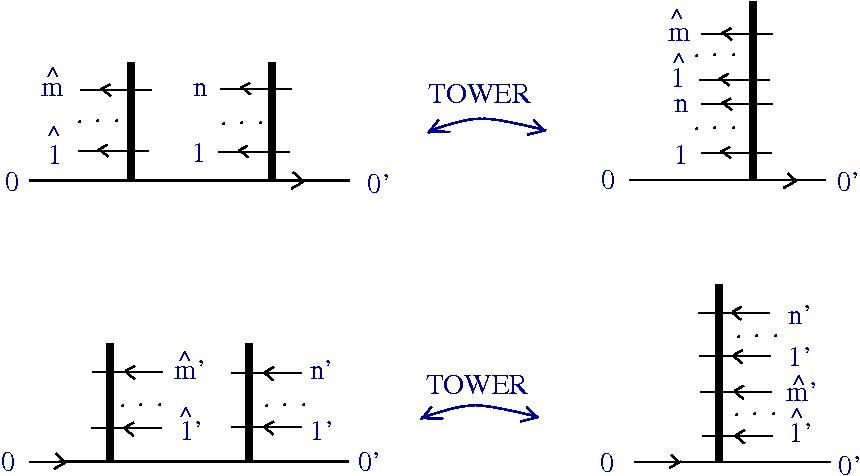}
     \caption{The TOWER moves.}
     \label{zipper_not_2}
 
\end{figure}

\item[-] CO-COMM, CO-ASSOC and FAN-IN moves, Fig.~\ref{convention_3}, are the same as the ones from chemlambda \cite{buligachem} \cite{alife}, see also section \ref{schemlambda}; 

\begin{figure}[H]
      
     \includegraphics[width=  95mm]{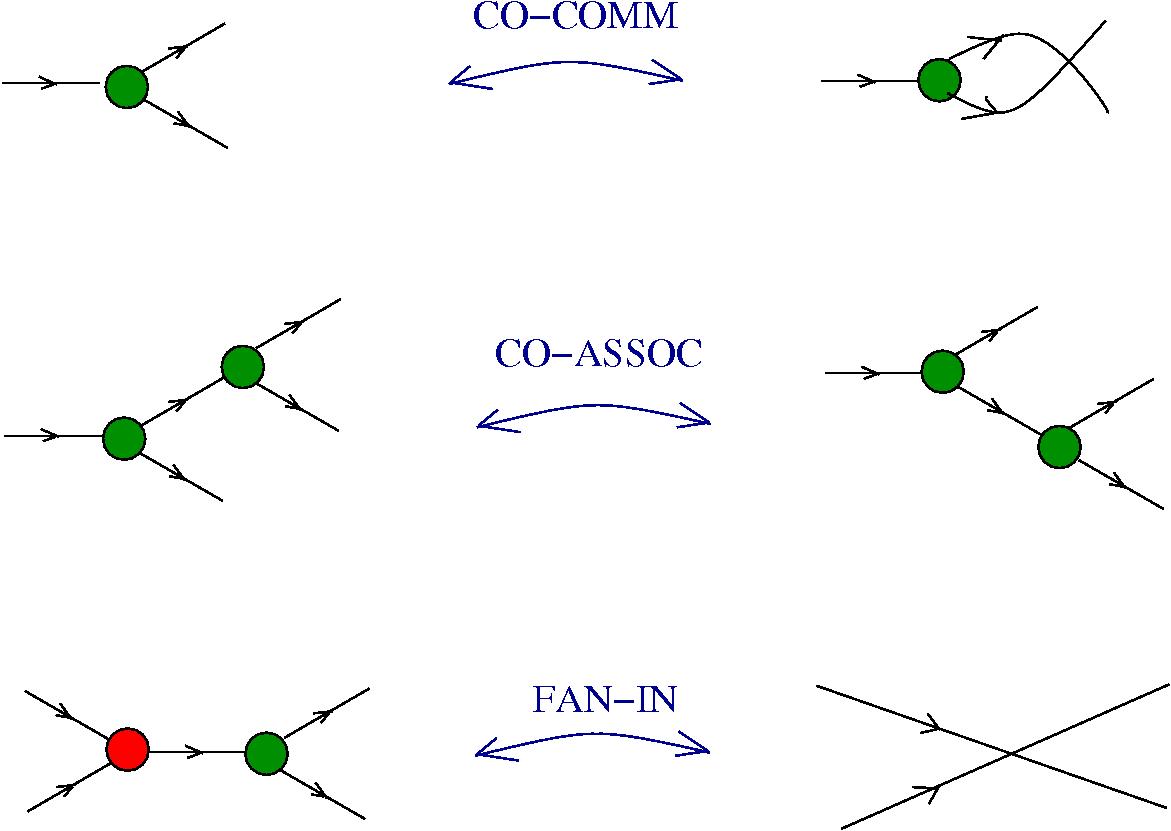}
     \caption{Some of the moves involving the other nodes: (1st row) CO-COMM move, (2nd row) CO-ASSOC move, (3rd row) FAN-IN move.}
     \label{convention_3}
 
\end{figure}

\item[-] the DIST moves are described in Fig.~\ref{zip_dist}.

\begin{figure}[H]
      
     \includegraphics[width=  110mm]{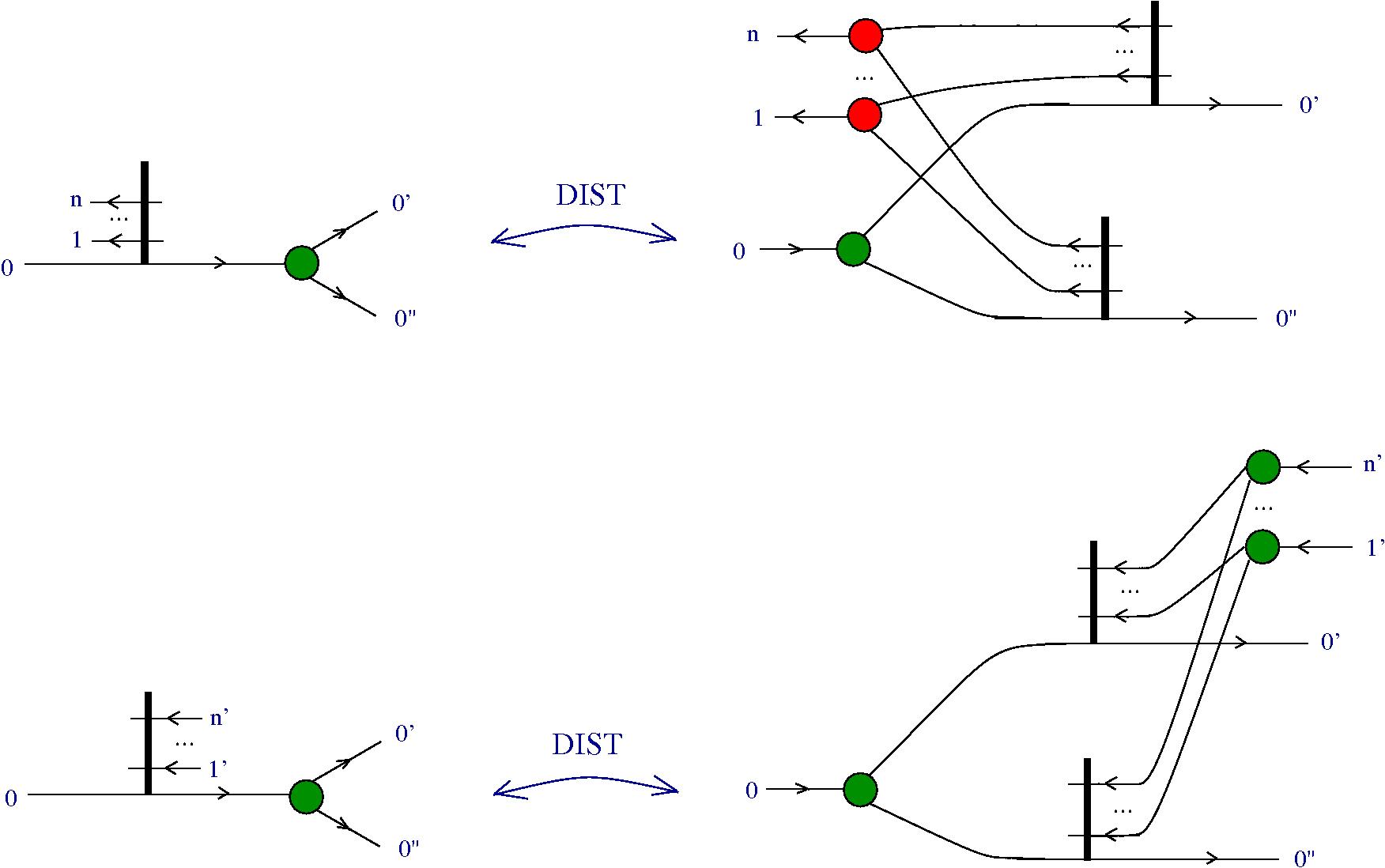}
     \caption{The DIST moves.}
     \label{zip_dist}
 
\end{figure}

\begin{figure}[H]
      
     \includegraphics[width=  95mm]{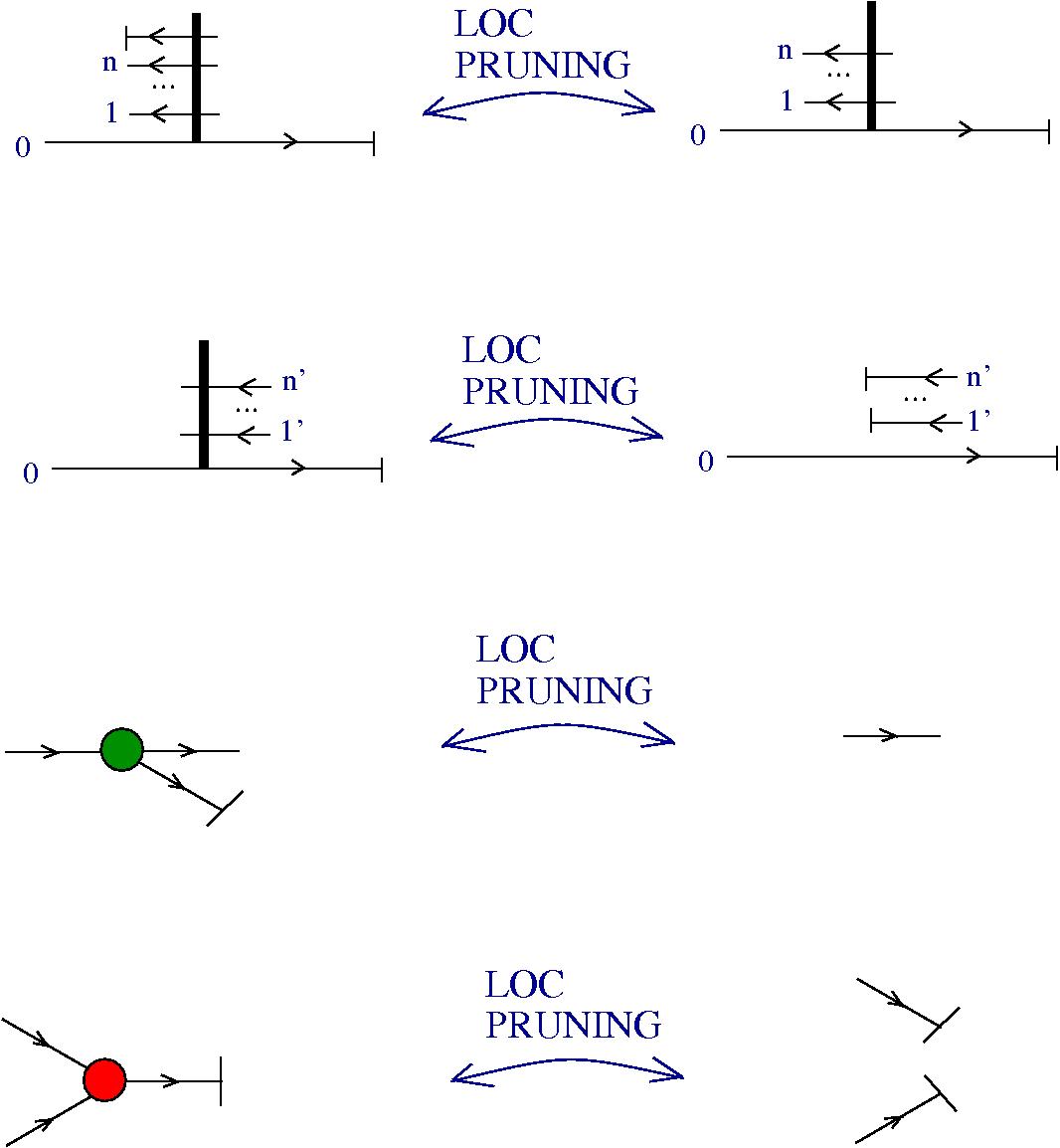}
     \caption{Local pruning moves. The first two rows describe moves for half-zippers, the last two rows describe moves for the other nodes.}
     \label{convention_4}
 
\end{figure}

\item[-] the last two rows of Fig.~\ref{convention_4} are the LOCAL PRUNING moves for the fanout and fanin nodes from chemlambda, the first two rows of the mentioned figure describe the LOCAL PRUNING moves of half-zippers;

\end{enumerate}
\label{defzippermoves}
\end{definition}

\section{Relations to chemlambda}
\label{schemlambda}

Chemlambda is a graph rewrite system, dressed as an artificial chemistry. It consists into a set of graphs called molecules and a collection of local moves, or graph rewrites. It has been introduced in \cite{buligachem} as an alternative to graphic lambda calculus, or GLC, \cite{bgraph}. In chemlambda there are only local moves, i.e. there is an a priori bound on the number of edges and nodes which are involved in any of the moves.

There is a distributed, decentralized computing model associated to chemlambda (or GLC), called distributed GLC 
\cite{distglc}. A good introduction to chemlambda, which also emphasizes it's biological like self-multiplication features, is \cite{alife}.

\begin{figure}[H]
      
     \includegraphics[width=  85mm]{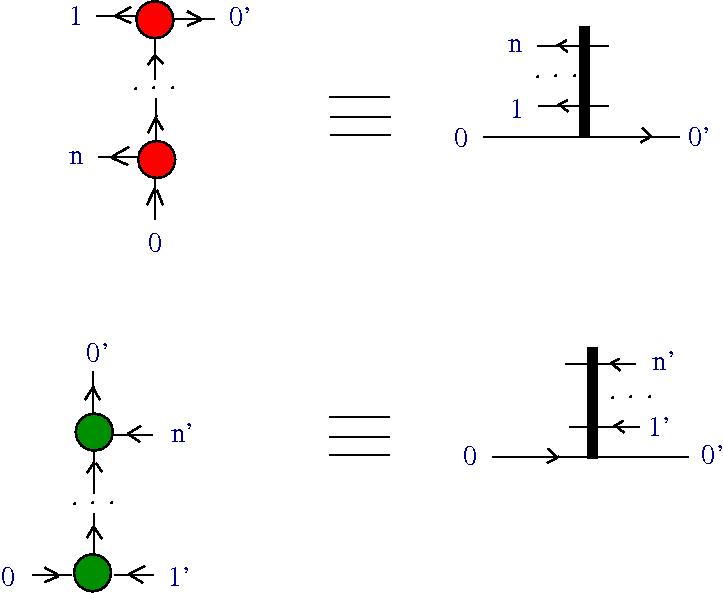}
     \caption{Zippers in chemlambda.}
     \label{zipper_not_1}
 
\end{figure}

\begin{proposition}
(a) In chemlambda, let's define half-zippers as in the Fig.~\ref{zipper_not_1} and let's use the CLICK move from the Fig.~\ref{click} as the definition of a zipper.  Then every move from the zipper logic can be realized as a finite sequence of chemlambda moves. In particular the DIST moves for half-zippers show that they are distributors, in a generalized sense, explained in \cite{alife} section "Propagators, distributors, multipliers and guns". 

(b) In zipper logic, let's define the lambda abstraction node as a $(-1)$ half-zipper, the application node as a $(+1)$ half-zipper. Then the moves  CLICK followed by ZIP is the beta move from chemlambda, the TOWER and CLICK moves serve to translate from zipper graphs to chemlambda molecules and the rest of the moves, used only with 1 half zippers,  are exactly the moves of chemlambda. 

\label{ptranslate}
\end{proposition}

We leave the proof to the reader, instead of giving detailed explanations,  because it is simply a matter of comparison of graph rewrites from the two formalisms. Please use \cite{alife} section "The Chemlambda formalism", Figures 1-5, as reference for the chemlambda moves.

\section{Zipper combinators}
\label{sturing}

From Proposition \ref{ptranslate} we get that zipper logic and chemlambda are equivalent. In particular, it follows that zipper logic is Turing universal, as chemlambda, because it contains combinatory logic. 

 However, the zipper logic may be more intuitive than chemlambda. Let's see how the SKI system of combinators appear and function in zipper logic. 

\begin{figure}[H]
      
     \includegraphics[width=  85mm]{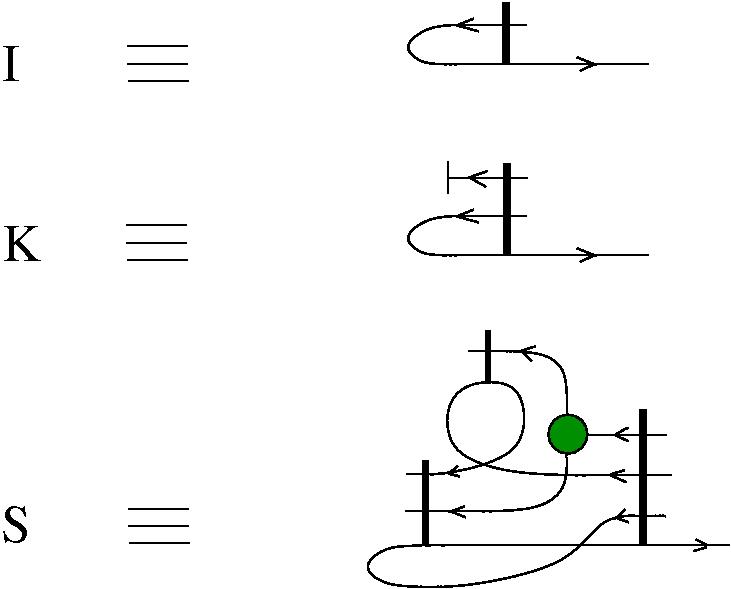}
     \caption{The S, K, I zipper combinators.}
     \label{zipper_not_6}
 
\end{figure}

\begin{definition}
The set of zipper combinators is the smallest set of zipper graphs with the properties: 
\begin{enumerate}
\item[-] it contains the S, K, I zipper graphs defined in Fig~\ref{zipper_not_6},

\item[-] for any natural number $n >0$ and for any $n+1$ zipper combinators, the zipper graph obtained by connecting the out arrows of the  zipper combinators to the in arrows of the $(+n)$ half-zipper is a zipper combinator.  
\item[-] any zipper graph which is obtained by applying a zipper logic move to a zipper combinator is a zipper combinator. 
\end{enumerate}
\label{dzippercombi}
\end{definition}

\begin{figure}[H]
      
     \includegraphics[width=  95mm]{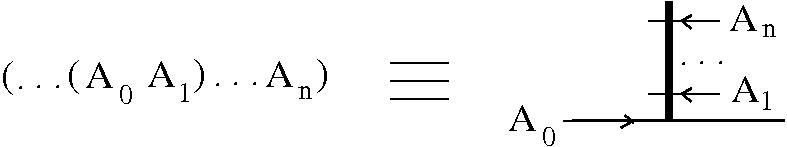}
     \caption{(+) half-zippers as concatenations of applications at right. This is compatible with Fig.~\ref{zipper_not_1}.}
     \label{appli}
 
\end{figure}

For any two zipper combinators $A$, $B$, the zipper combinator $AB$ is obtained by connecting $A$ and $B$ to a $(+1)$ half-zipper. More generally, from any  $n+1$ zipper combinators $A_{0}$, ... , $A_{n}$, we obtain the zipper combinator $(...(A_{0}A_{1})...A_{n})$ described in Fig.~\ref{appli}. 

\paragraph{Birth and death of zipper combinators.} We shall need the following lemmata, which describe how zipper combinators multiply or die.

\begin{lemma}
Any zipper combinator is a multiplier, i.e. for any zipper combinator $A$, the graph obtained by connecting the out arrow of $A$ to the entry arrow of a fanout node transforms by a finite sequence of zipper logic moves into two copies of $A$. 
\label{lmulti}
\end{lemma}

\paragraph{Proof.} Suppose that the zipper combinator $A$ is formed by    $n+1$ zipper combinators $A_{0}$, ... , $A_{n}$ connected to a $(+n)$ half-zipper. Use then the DIST move for $(+)$ half-zippers from Fig.~\ref{zip_dist}, and remark that if $A_{0}$, ... , $A_{n}$ are multipliers then $A$ is a multiplier.  From the definition \ref{dzippercombi} it follows that in order to prove our lemma, it is sufficient to prov that $S$, $K$ and $I$ zipper combinators are multipliers. But this has been shown several times in previous articles. This is related to  \cite{alife} section "Propagators, distributors, multipliers and guns". In the Step 2 of the  proof of Theorem 4.2 
\cite{buligachem} there is such a detailed proof, only the formalisms differ slightly, i.e. in the mentioned reference it is used chemlambda, while here we use zipper logic, moreover there was used the BCKW system of combinators, while here it is used the SKI system of combinators. That is why we leave the completion of the rest of the proof to the interested reader. \hspace{.5cm} $\square$

\begin{figure}[H]
      
     \includegraphics[width=  85mm]{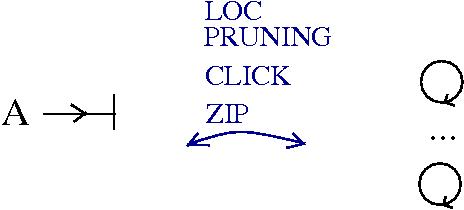}
     \caption{Every zipper combinator connected to a termination node transforms into a finite collection of loops.}
     \label{zipd_5}
 
\end{figure}

\begin{lemma}
For any zipper combinator $A$, the graph obtained by connecting the out arrow of  $A$ with a termination node  can be reduced, by a finite number of moves, to a graph formed by a finite number of loops (or the empty graph). See Fig.~\ref{zipd_5}. 
\label{ldie}
\end{lemma}

\paragraph{Proof.} Consider the zipper combinator obtained from  $n+1$ zipper combinators $A_{0}$, ... , $A_{n}$ connected to a $(+n)$ half-zipper, then connect the out arrow of this zipper combinator to a termination node. By the LOC PRUNING move from Fig.~\ref{convention_4}, second row,  we can reduce this graph to the collection of  $A_{0}$, ... , $A_{n}$, each connected to a termination node.

\begin{figure}[H]
      
     \includegraphics[width=  95mm]{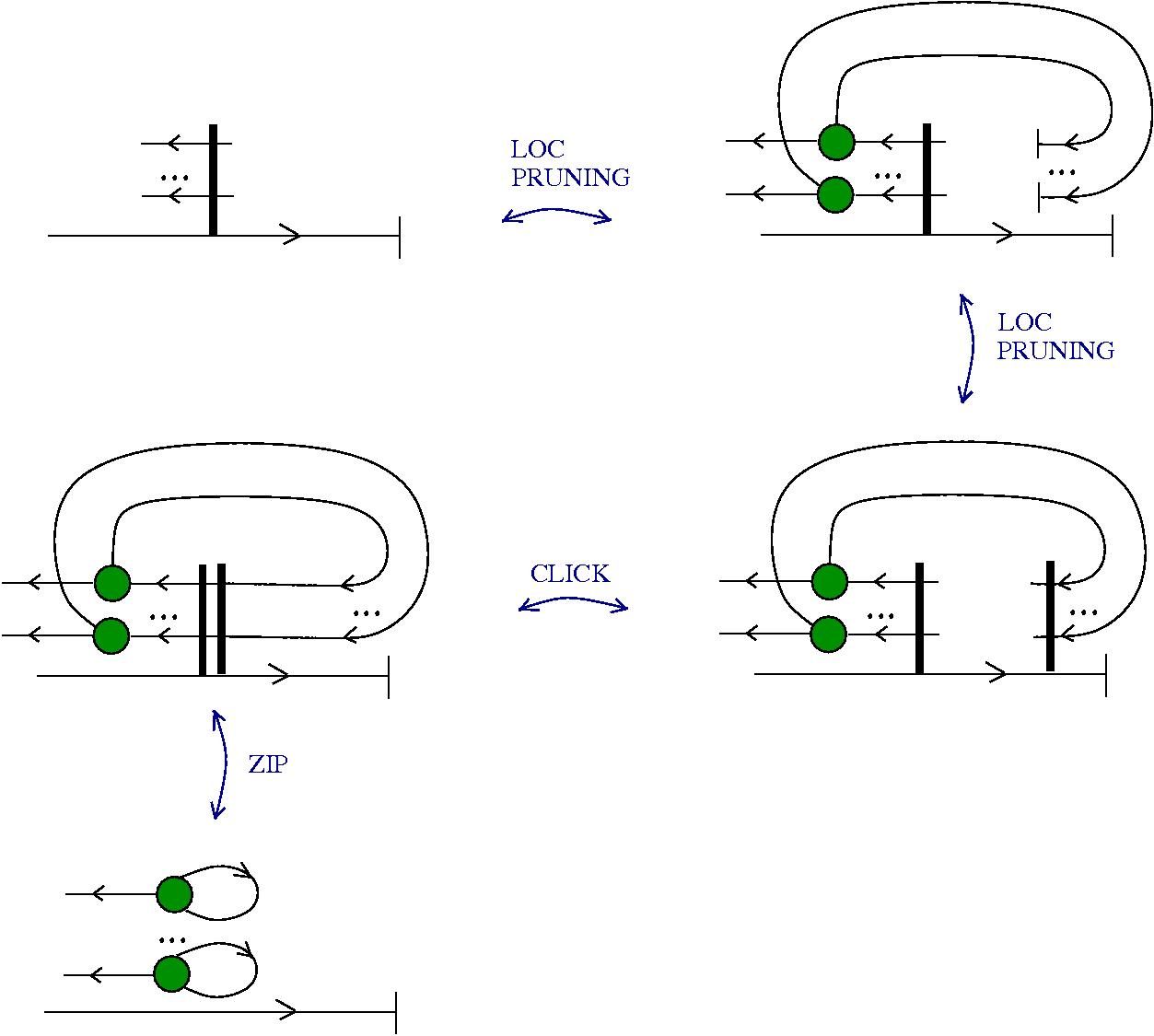}
     \caption{A succession of moves used to eliminate a (+) half-zipper connected to a termination node.}
     \label{zipd_6}
 
\end{figure}

\begin{figure}[H]
      
     \includegraphics[width=  95mm]{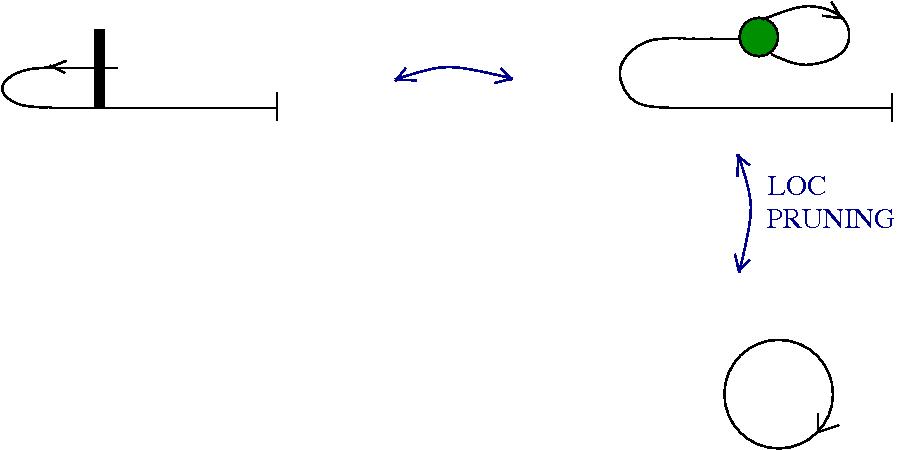}
     \caption{The $I$ combinator connected to a termination node transforms into one loop.}
     \label{zipd_7}
 
\end{figure}

For progressing further we shall use the following "trick", described in Fig.~\ref{zipd_6}, which will be used to pass over a (+) half-zipper connected to a termination node.

\begin{figure}[H]
      
     \includegraphics[width=  105mm]{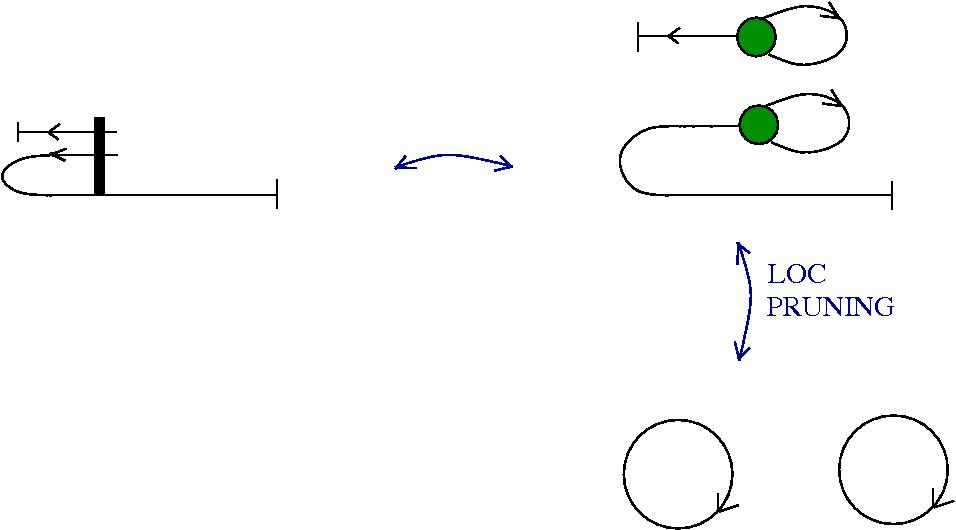}
     \caption{The $K$ combinator connected to a termination node transforms into two loops.}
     \label{zipd_8}
 
\end{figure}

\begin{figure}[H]
      
     \includegraphics[width=  125mm]{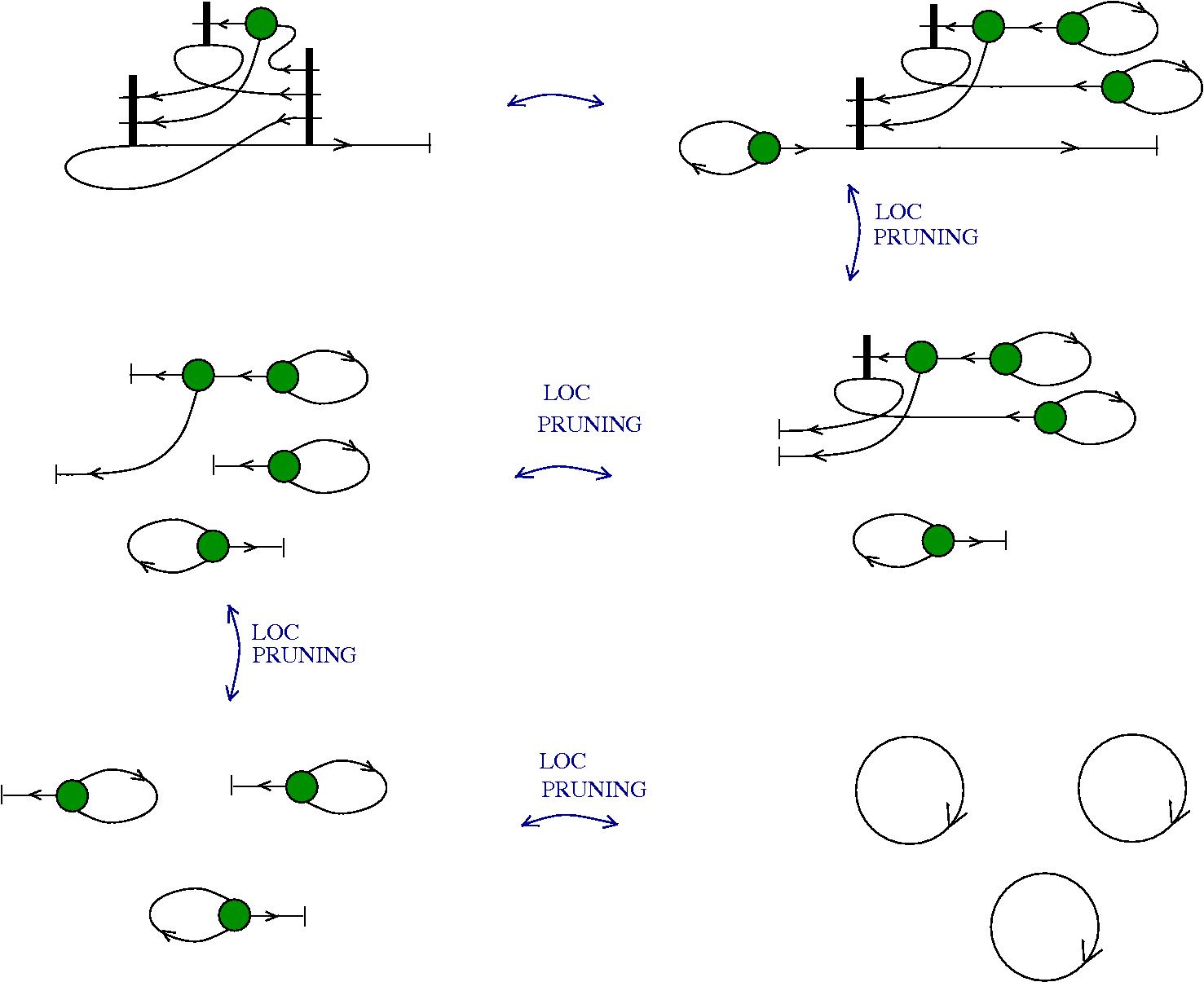}
     \caption{The $S$ combinator connected to a termination node transforms into three loops.}
     \label{zipd_9}
 
\end{figure}

In the second part of the proof we use the "trick" for the zipper combinators $I, K$ and $S$ connected to a termination node. In the case of the $I$ combinator. we can transform it into one loop, as described in Fig.~\ref{zipd_7}. The first blue arrow represents the "trick". 

The same trick is used for transforming the $K$ zipper combinator, connected to a termination node, into two loops, Fig.~\ref{zipd_8}.

The zipper combinator $S$, connected to a termination node, is transformed into three loops, starting by the same trick, Fig.~\ref{zipd_9}. \hspace{.5cm} $\square$

\begin{theorem}
Consider the set of all zipper combinators, with the relation "$A=B$" meaning that there is a finite sequence of zipper logic moves which transforms $A$ into $B$ and a finite number of loops, and with the operation $AB$ as defined in the Fig.~\ref{appli}. Then we have the relations: 
\begin{enumerate}
\item[(a)] $IA = A$ for any zipper combinator $A$, 
\item[(b)] $(KA)B = A$ for any zipper combinators $A, B$, 
\item[(c)] $(((SA)B)C) = (AC)(BC)$ for any zipper combinators $A, B, C$, 
\item[(d)] $(SK)K = I$. 
\end{enumerate}
\label{tzipc}
\end{theorem}

\paragraph{Proof.} We start by proving (a), (b) and (d). Then we prove (c) as a consequence of the fact that zipper combinators are multipliers. 
 
\begin{figure}[H]
      
     \includegraphics[width=  75mm]{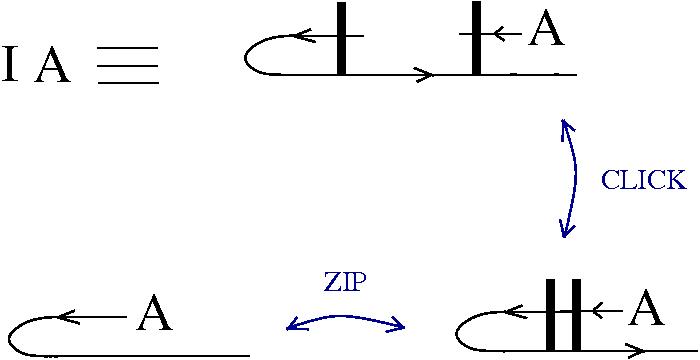}
     \caption{Proof of $IA = A$.}
     \label{zipper_not_7}
 
\end{figure}

The proof of (a) is given in Fig.~\ref{zipper_not_7}. 

\begin{figure}[H]
      
     \includegraphics[width=  95mm]{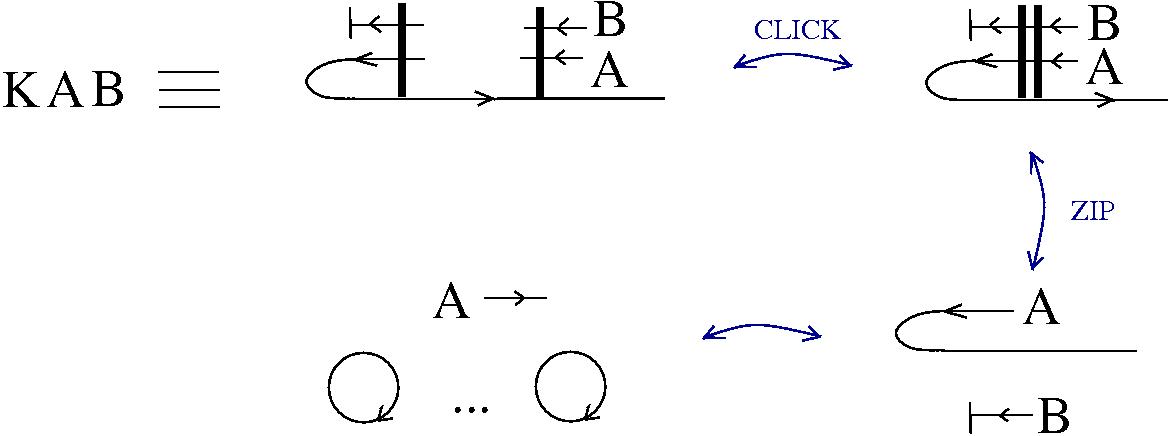}
     \caption{Proof of $(KA)B = A$.}
     \label{zipper_not_10}
 
\end{figure}

The proof of (b) is given in Fig.~\ref{zipper_not_10}. On the second row of the figure, the move from right to left is the one from Lemma \ref{ldie}, which transforms the  zipper combinator $B$ connected to a termination node into a finite collection of loops. 

\begin{figure}[H]
      
     \includegraphics[width=  130mm]{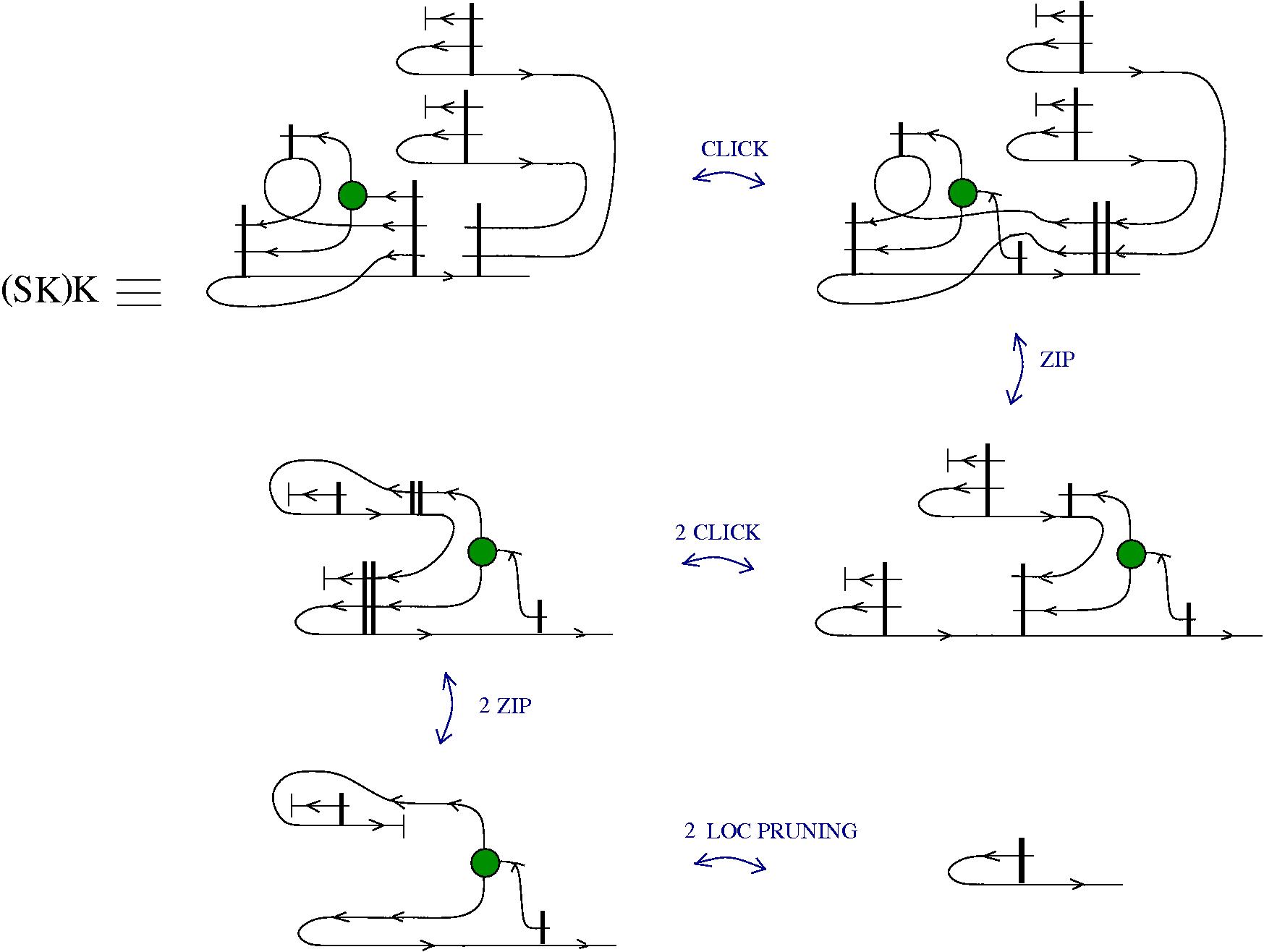}
     \caption{Proof of $(SK)K = I$.}
     \label{zipper_not_11}
 
\end{figure}

The proof of (d) is given in Fig.~\ref{zipper_not_11}.

\begin{figure}[H]
      
     \includegraphics[width=  125mm]{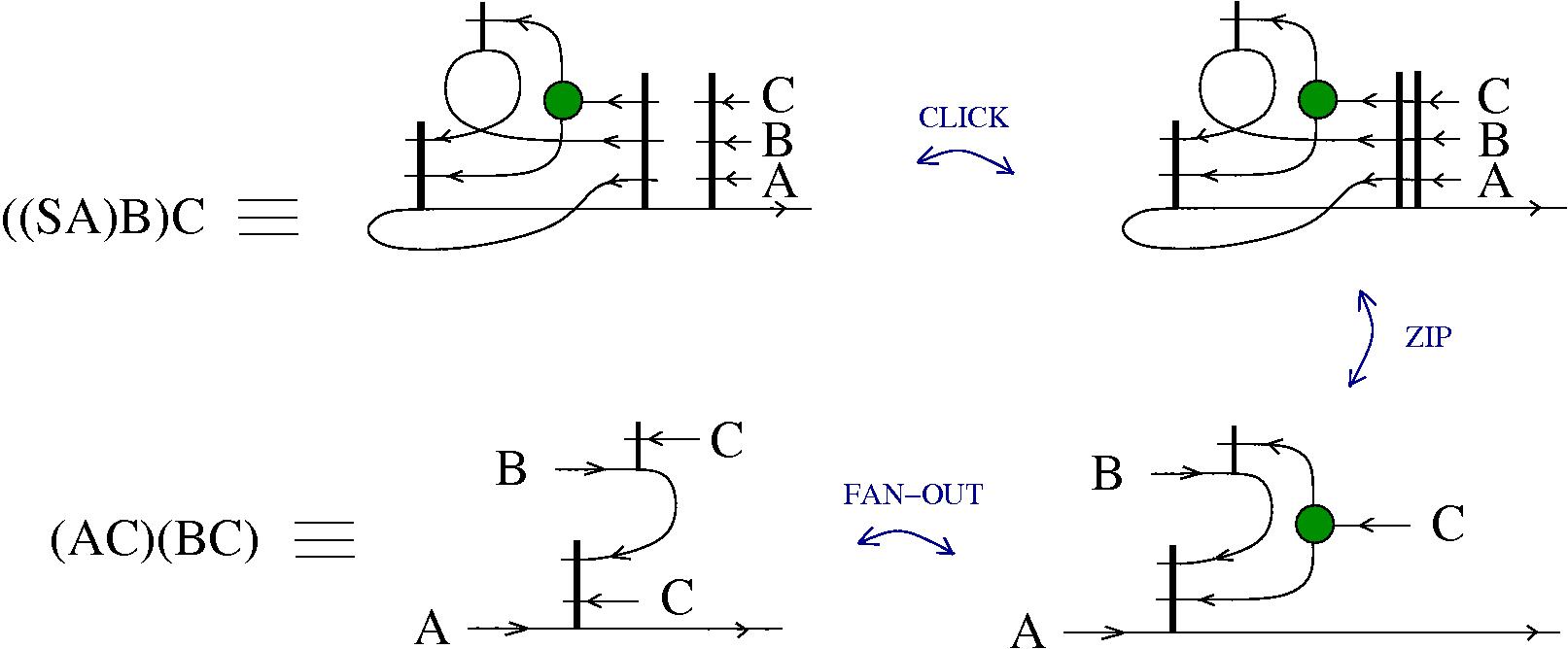}
     \caption{Proof of $((SA)B)C = (AC)(BC)$.}
     \label{zipper_not_12}
 
\end{figure}

For the proof of (c) see Fig.~\ref{zipper_not_12}.  

The "FAN-OUT" move from the second row of the figure is the one which multiplies the zipper combinator $C$, according to the lemma \ref{lmulti}. \hspace{.5cm} $\square$

\section{Zippers and knots}
\label{sknots}

The fact that there are relations between knots and tangles diagrams and graph rewrites systems which are close to zipper logic has been already noticed. In the article which introduces GLC (or "graphic lambda calculus") \cite{bgraph}, in the section 6, see also Theorem 6.1, is proved that GLC has a sector which is equivalent with the formalism of locally planar oriented crossings diagrams with the oriented Reidemeister moves. In this sector we represent an oriented crossing as a pair of two nodes, lambda abstraction and application. In \cite{distglc}, section 5, are discussed various connections between a topological version of GLC and knot diagrammatic notations for lambda calculus or for topological quantum computations. In the last section of \cite{alife}, there is proposed another encoding of an oriented crossing, as a pair fanout and application nodes, which leads to a graph rewrite formalism over oriented knots or tangle diagrams, with a lambda abstraction trivalent node added, along with a termination univalent node.

\begin{figure}[H]
      
     \includegraphics[width=  115mm]{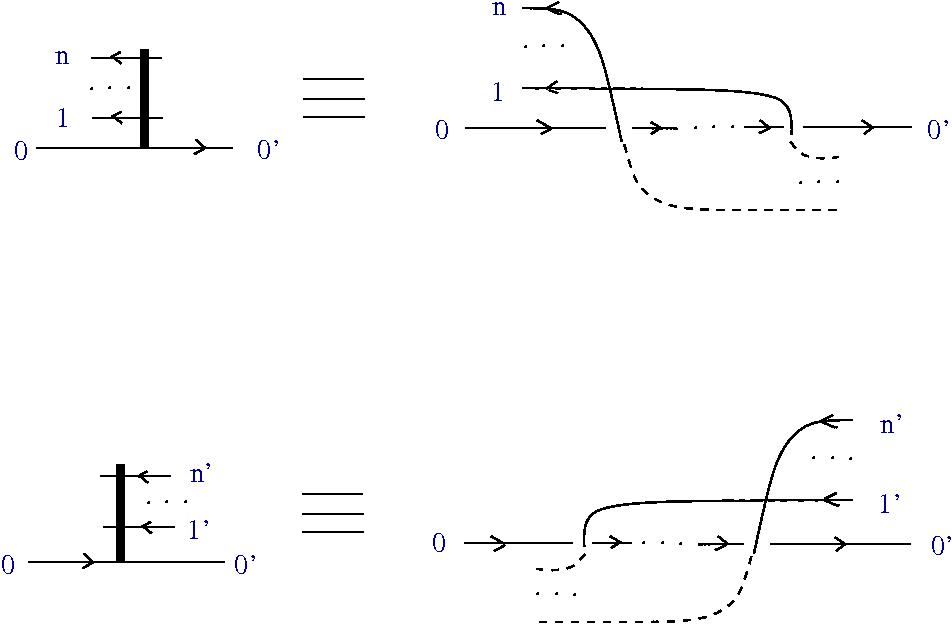}
     \caption{Definition of zippers as knot diagrams.}
     \label{zipper_loop_1}
 
\end{figure}

The zipper logic formalism can be transformed into a graph rewriting formalism over knots or tangles diagrams, with the trivalent fanin and fanout nodes added, along with the univalent termination node. Indeed, let's define half-zippers like in the Fig.~\ref{zipper_loop_1}.

\begin{figure}[H]
      
     \includegraphics[width=  95mm]{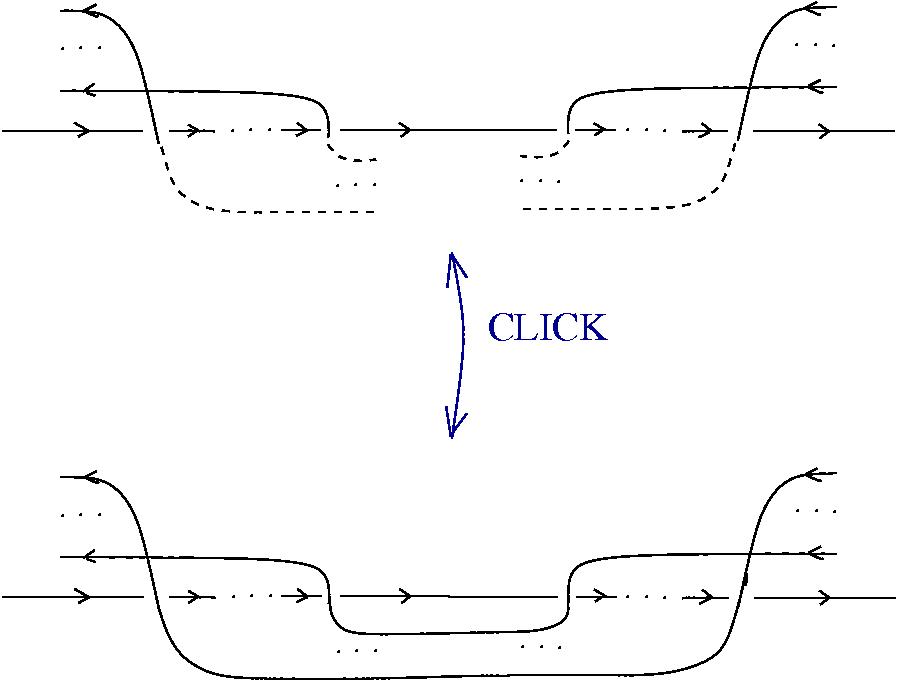}
     \caption{The CLICK move for zippers as knot diagrams.}
     \label{zipper_loop_2_n}
 
\end{figure}

The dotted lines represent arcs which are only virtually there. Their advantage is that we can see the CLICK move, in the realm of knot diagrams, as a move which transforms virtual arcs into real, connected arcs, see Fig.~\ref{zipper_loop_2_n}. 

This way, the ZIP move (which is the equivalent of the graphic beta move from GLC) appears as one of the oriented Reidemeister 2 moves.

\begin{figure}[H]
      
     \includegraphics[width=  100mm]{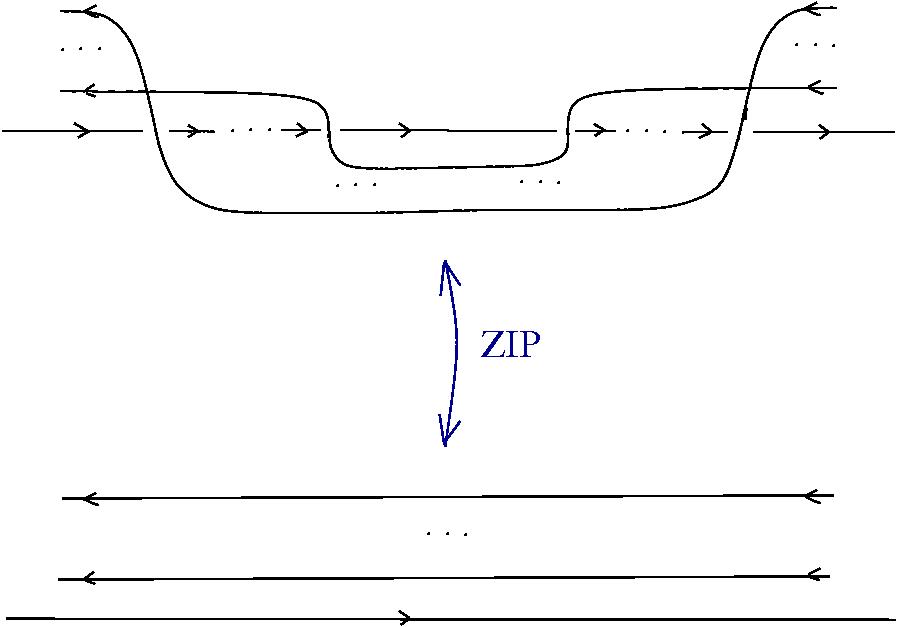}
     \caption{The ZIP move for zippers as knot diagrams is the same as R2.}
     \label{zipper_loop_3_n}
 
\end{figure}

\begin{figure}[H]
      
     \includegraphics[width=  80mm]{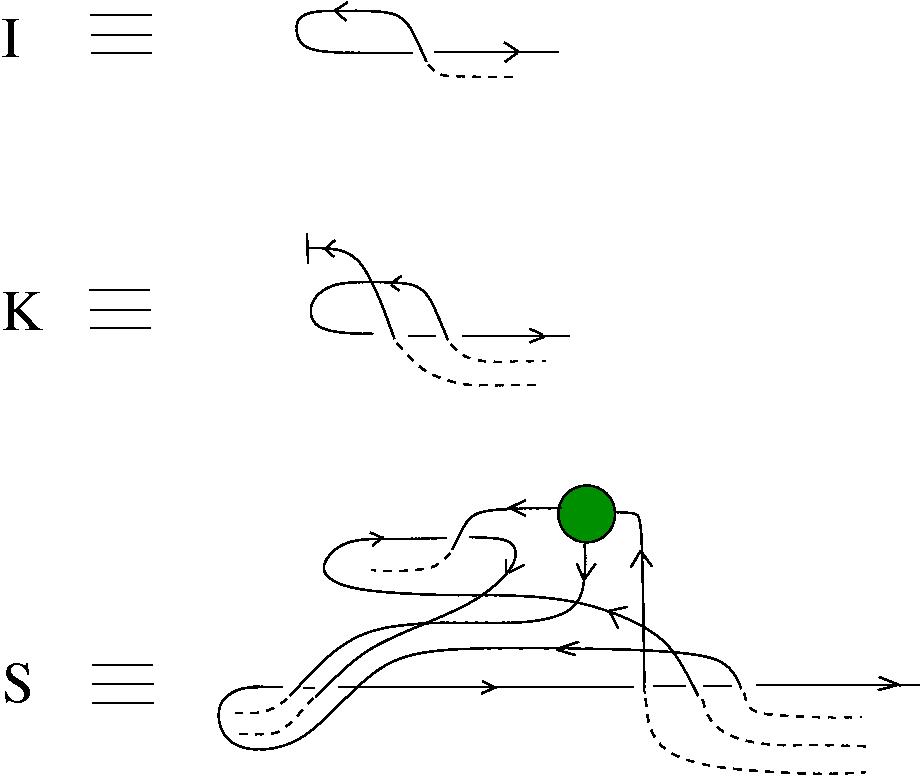}
     \caption{The S,K,I combinators as knot diagrams.}
     \label{zipper_loop_4}
 
\end{figure}

Finally, the S,K,I zipper combinators appear as the following "knot combinators" from Fig.~\ref{zipper_loop_4}. 

The connections between knot diagrammatics and lambda calculus started by the article Kauffman \cite{kauffman2}, where in section 5 knot diagrams are used as notations for combinatory logic terms. See also \cite{kauffman3} where knot automata are introduced. The most recent discussion concerning diagrammatic methods for representing quantum processes and quantum computing is to be found in \cite{kauffman5}. Another interesting, related research thread is the one of tangle machines by Carmi and Moskovich \cite{carmimosko1} \cite{carmimosko2}.  In all these articles knot diagrams are used as a notational device for computations. The topology appears as related to the invariance of these computations, in the sense that knot diagrams which are related by a sequence of Reidemeister moves describe the same computation. In few words, in this case we may say that "topology does not compute, but is an invariant". 

In contradistinction, in the knot diagrams sector of GLC, or in the topological version of GLC the topology does compute. This means that (some of the) reduction moves appear as graph rewrites on knot diagrams which change their topology. For example the graphic beta move is a crossing smoothing move in the knot diagrams sector of GLC. 

It would be interesting to explore in more detail this distinction -- topology does not compute vs topology does compute -- for the benefit of all these diagrammatic formalism. The distinction may turn out to be more subtle. For example, Kauffman' bracket polynomial algorithm \cite{kauffman4} uses  a combination of skein relations and smoothing moves, thus, as remarked in \cite{distglc}, section 5, "This is similar to allowing free beta reduction in the lambda calculus graphs. [...] An analogous situation could occur in GLC   where one would need the average over all the results of
the many branching calculations."  There exist probably overarching formalisms, which blend the two roles of topology in knot diagrammatics related to logic, waiting to be discovered.

\end{document}